\documentclass[%
preprint,
preprintnumbers,
 amsmath,
 amssymb,
]{revtex4-2}

\usepackage{graphicx}
\usepackage{dcolumn}
\usepackage{bm}
\usepackage[mathlines]{lineno}

\begin{document}

\title{Gain modulation unlocks phase plasticity, thus allowing for malleable nonlinear dynamics}

\author{Christoforos A. Papasavvas}
\email{christoforos.papasavvas@ncl.ac.uk}

\affiliation{%
 School of Computing, Newcastle University, UK
}%

\begin{abstract}

Nonlinear dynamics emerge through either nonlinear interactions between the variables or through nonlinearities imposed on their linear interactions. Their interactions can be conceptualized as modulations of input-output (I/O) functions, where one variable modulates another variable’s I/O function. Within this framework, nonlinear interactions are manifested as gain modulations, where the gain (i.e. slope) of a variable’s I/O function is modulated by another variable. By analyzing oscillatory dynamics in modulation phase space, I show that gain modulation qualitatively enhances phase plasticity by loosening the relationship between oscillation phase and instantaneous gain. This finding reveals how gain modulation renders nonlinear dynamics more adaptable. I discuss the direct implications of this finding on the ability of any nonlinear system to get entrained to rhythmic input and synchronize with other systems and how it relates to neural dynamics.


\end{abstract}

\date{\today}

\maketitle


Multivariable dynamical systems are defined by sets of differential equations which express the time evolution of the state variables. Even a single nonlinear term in those equations renders the system nonlinear. Nonlinear terms appear in two distinct forms. First, they can appear as nonlinear functions imposed on single variables or on a linear combination of multiple variables (i.e., $\bm{u}(p_1v_1+p_2v_2+...)$, with parameters $p$, variables $v$, and nonlinear function $\bm{u}$). Second, they can appear as nonlinear interactions (e.g., multiplications) between the variables or their functions (e.g., $\bm{u_1}(v_1)\times \bm{u_2}(v_2)\times ...$, with variables $v$ and linear/nonlinear functions $\bm{u}$). In the latter, the variables interact in a nonlinear fashion which, as explained below, can be considered as \textit{gain modulation}; whereas in the former, any nonlinearity is imposed onto their linear interactions (conceptualized as \textit{linear modulations} below). For instance, nonlinear dynamical systems that feature nonlinearities exclusively in the form of imposed nonlinear functions include the FitzHugh–Nagumo model \cite{Fitzhugh1961} and the Goodwin oscillator \cite{Goodwin1965, Goldbeter1995}. In contrast, the R\"ossler system \cite{Rossler1976} and the predatory-prey system \cite{Freedman1980} feature nonlinear interactions (i.e., gain modulation) between the variables.

Typical nonlinear phenomena, such as limit cycles and chaos, can be produced through either of those types of nonlinearities. Chaos, for instance, emerges through the nonlinear interaction between the variables in the R\"ossler system \cite{Rossler1976}, while it relies on univariate nonlinear terms in the Hindmarsh–Rose system \cite{Hindmarsh1984, Osipov1998}. Limit cycles are also found in both of those systems, without any obvious qualitative difference between them. However, a recent finding suggests that the adaptability and entrainability of such limit cycles depend on the type of nonlinearities featured in the system, with gain modulation enhancing entrainment \cite{Papasavvas2020}. Furthermore, the entrainability of limit cycles has been reported to depend on their adaptability and phase plasticity, that is, their ability to phase shift in order to accommodate parameter changes and transient stimuli \cite{Hatakeyama2015, Paijmans2017}. Motivated by these recent findings, the present Letter addresses the following question: Is the adaptability of nonlinear dynamics dictated by the type of nonlinearities featured in the system?

\section*{Conceptualization}

In order to investigate this, I consider the original R\"ossler system and two variations of it. The R\"ossler system, Eqs.~(\ref{eq1}), is a well studied nonlinear 3-dimensional system producing limit cycles, bifurcations, and chaos \cite{Rossler1976}. Despite its rich behavior, it is a particularly elegant system, since there is only one nonlinear term, $\bm{f}(x,z)$. Originally, $\bm{f_0}(x,z) = xz$. The two variations of the system considered here are defined with an imposed nonlinearity (justification in Generalization). A sigmoidal function, specifically, the hyperbolic tangent, is imposed onto two distinct types of interactions between the variables $x$ and $z$:  $\bm{f_1}(x,z) = \tanh(x+z)$ and $\bm{f_2}(x,z) = \tanh(xz)$. In the first variation, the variables interact linearly, whereas in the second they interact nonlinearly. 

\begin{subequations}
\label{eq1}
\begin{equation}
\frac{\mathrm{d}x}{\mathrm{d}t}= -y -z
\label{subeq:1}
\end{equation}
\begin{equation}
\frac{\mathrm{d}y}{\mathrm{d}t}= ax+by
\label{subeq:2}
\end{equation}
\begin{equation}
\frac{\mathrm{d}z}{\mathrm{d}t}= c -dz +\bm{f}(x,z)
\label{subeq:3}
\end{equation}
\end{subequations}

The traditional way of interpreting the first variation of the nonlinear term, $\bm{f_1}(x,z) = \tanh(x+z)$, is the following: the sum of variables $x$ and $z$ serves as the input to the hyperbolic tangent. That is, $\bm{f_1}(x,z) = \bm{g}(\bm{h_1}(x,z))$, where $\bm{g}(\cdot)=\tanh(\cdot)$ and $\bm{h_1}(x,z)=x+z$. An alternative interpretation considers $\bm{g}$ as the input-output (I/O) function of one of the variables, while the other variable modulates it \cite{Silver2010}. Here, the I/O function $\bm{g}$ of $z$ is linearly modulated by $x$, or vice versa. That is, $\bm{f_1}(x,z) = \bm{g}(z+\bm{h}(x))$, where $\bm{h}(x)=x$ and $\bm{g}$ as above. Here, $\bm{h}$ serves as the linear modulator of $\bm{g}$'s input. This modulation is expressed by shifting the whole sigmoidal $\bm{g}$ along the input axis. 

Similarly, the traditional way of interpreting the second variation of the nonlinear term, $\bm{f_2}(x,z) = \tanh(xz)$, is as follows: the product of variables $x$ and $z$ serves as the input to the hyperbolic tangent. That is, $\bm{f_2}(x,z) = \bm{g}(\bm{h_2}(x,z))$ where $\bm{h_2}(x,z)=xz$ and $\bm{g}$ as above. The alternative interpretation is: the I/O function $\bm{g}$ of $z$ is gain modulated by $x$, or vice versa \cite{Silver2010}. That is, $\bm{f_2}(x,z) = \bm{g}(z\bm{h}(x))$, with $\bm{h}$ and $\bm{g}$ as above. Now, $\bm{h}$ serves as the gain (nonlinear) modulator of $\bm{g}$'s input. This nonlinear modulation is expressed by squeezing and stretching the I/O function, thus inevitably changing its slope (i.e. gain). This is precisely how gain modulation is defined here: any modulation of the function's slope \cite{Silver2010}. 

I will use the alternative interpretations throughout this Letter, thus nonlinear dynamics are considered to be shaped through the modulation of $z$'s I/O function by $x$.  Therefore, the first variation of the system is characterized by a linear modulation of a nonlinear I/O function due to the linear interaction between the variables, whereas the second variation features gain modulation of the same I/O function due to the nonlinear interaction between the variables.

\section*{Analysis in modulation phase space}

\begin{figure*}
\includegraphics{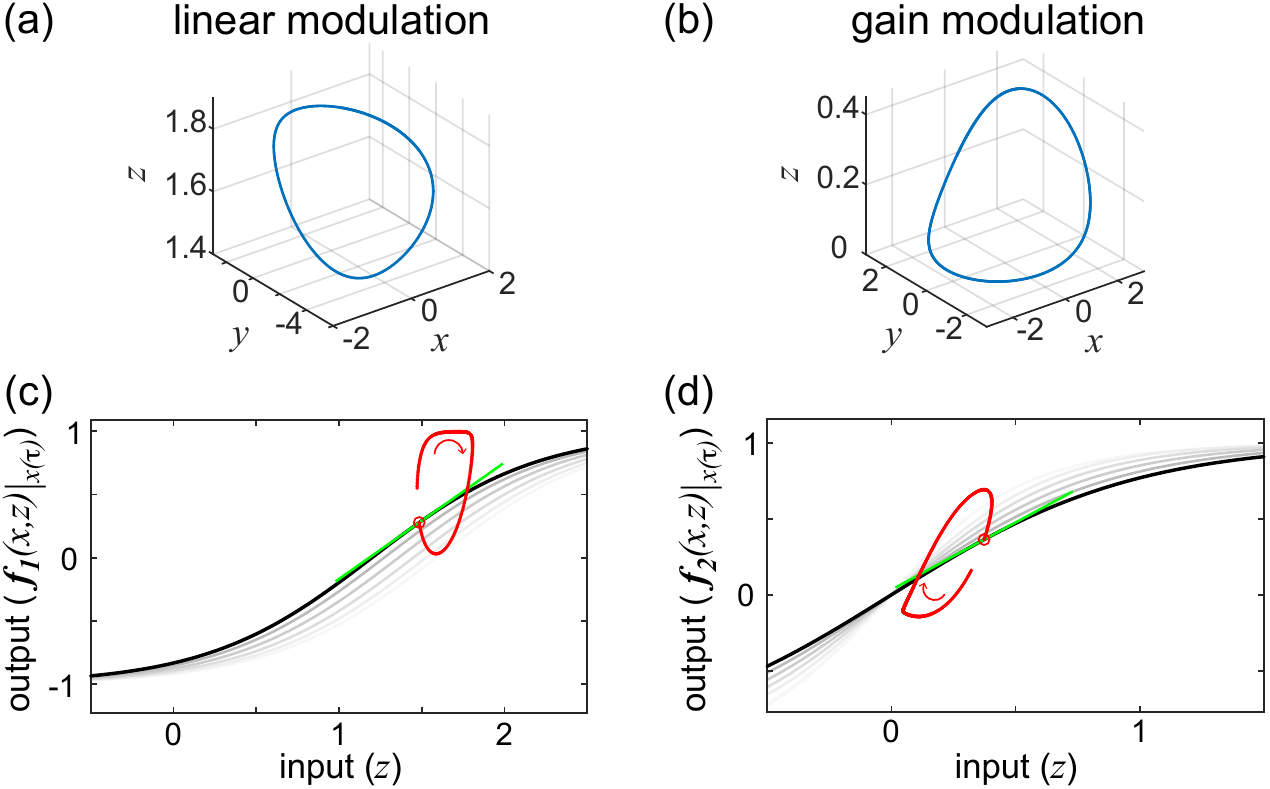}
\caption{\label{fig1}
Limit cycles in two variations of the R\"ossler system, one featuring linear modulation and the other gain modulation of the same sigmoidal I/O function. (a) A limit cycle in the $x$-$y$-$z$ phase space of the first variation. (b) Respective presentation of a limit cycle for the second variation. (c) Projection of the limit cycle (red curve) on the input-output space for the first variation. Notice the linear modulation of the sigmoidal, indicated by shaded black curves. The limit cycle is formed through this modulation while the state of the system at time $\tau$ is indicated by the red circle. The sigmoidal's tangent at this specific state is shown in green. (d) Respective projection for the second variation. Notice the mechanism of gain modulation, indicated by shaded black curves, shaping the limit cycle (red curve). See also animations in Supplementary Information.
}
\end{figure*}

Now, let us analyze stable oscillatory behavior in these two variations of the system by dissecting their limit cycles in phase space. Since limit cycles are nonlinear phenomena, the nonlinear term $\bm{f}(x,z)$ is the only term that can be responsible for the emergence of such phenomena in the system. Thus, the analysis will be focused on this term. For the first variation, let us consider the limit cycle produced with parameters $a=1$, $b=0.06$, $c=0.2$, and $d=2$. For the second variation, parameters $a=6$, $b=0.03$, $c=1.1$, and $d=1.1$ were used. Figure~\ref{fig1}(a)-(b) shows the two limit cycles in the $x$-$y$-$z$ phase space. The limit cycles seem to be qualitatively indistinguishable. What follows in this Letter is a demonstration that these limit cycles have fundamentally different quality due to the fact that only one of them emerges through gain modulation.

By considering the two different modulations applied to the sigmoidal I/O function, we can examine the formation of both limit cycles in the input-output space of the nonlinear term $\bm{f}(x,z)$, as shown in Fig.~\ref{fig1} (c)-(d). The input to the sigmoidal is $z$ alone, whereas variable $x$ serves as the modulator, either linear (Fig.~\ref{fig1}(c)) or nonlinear (Fig.~\ref{fig1}(d)). The plotted position and shape of the sigmoidal is evaluated for a specific modulator value at time $\tau$, $\bm{f}(x,z)\Bigr\rvert_{x(\tau)}$. The shading around the sigmoidal indicates the corresponding modulation in each case, shifting vs squeezing. Separating the two variables into input and modulator reveals the limit cycles in I/O space, formed by the red curve which traces the evolution of the system. It is important to consider the different modulations of the I/O function here, following the alternative interpretation of $\bm{f}(x,z)$ (see above). Mathematically and numerically the system behaves the same way as in the traditional interpretation. Yet, conceptually it is different and it allows for the graphical decomposition of the dynamics at the level of individual I/O functions.

Let us now consider the tangent of the sigmoidal at each point along the limit cycle. The tangent is shown in green in Fig.~\ref{fig1} (c)-(d). I define the instantaneous gain, $G$, as the slope of this tangent. More precisely, it is the partial derivative of the I/O function with respect to its input variable, in this case $G = \frac{\partial \bm{f}(x,z)}{\partial z}$. It expresses the instantaneous sensitivity of the output to input changes. So now, the modulation phase space of the nonlinear term can be defined as the 3-dimensional space of input-output-instantaneous gain (I/O/G). 

\begin{figure*}
\includegraphics{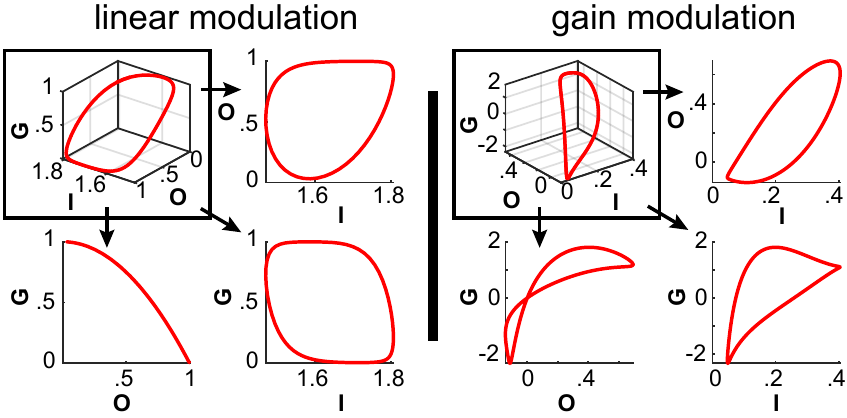}
\caption{\label{fig2}
Limit cycles in modulation phase space (I/O/G) for the two variations of the system. The modulation phase space is viewed from different angles in the subpanels, showing two dimensions at a time. Notice the qualitative difference between the two variations: limit cycles evolving on a 2-dimensional (left) vs 3-dimensional (right) manifold in modulation phase space.
}
\end{figure*}

The limit cycles in Fig.~\ref{fig1} are projected in the modulation phase space shown in Fig.~\ref{fig2}. By viewing the limit cycles from different angles, taking two dimensions at a time (subpanels in Fig.~\ref{fig2}), we observe a qualitative difference between them. The one produced through linear modulation has the dimensions of $G$ and O collapsed onto each other. More specifically, there is an one-to-one mapping from output to instantaneous gain, with a quadratic relationship between them, since $G = \frac{\partial (\tanh (x+z))}{\partial z} = \mathrm{sech}^2 (x+z) = 1 - \bm{f_1}^2(x,z)$. This relationship restricts the limit cycle on a 2-dimensional manifold in modulation phase space (see below for generalization). In contrast, the limit cycle produced through gain modulation evolves on a 3-dimensional manifold in modulation phase space.

At this point, I would like to highlight the interpretation of this qualitative difference in respect to the phase of the oscillation. For the limit cycle produced through linear modulation, a specific output always corresponds to the same instantaneous gain, regardless of the oscillation's phase. In contrast, for the limit cycle produced through gain modulation, the relationship between output and instantaneous gain always depends on the oscillation's phase.

\begin{figure*}
\includegraphics{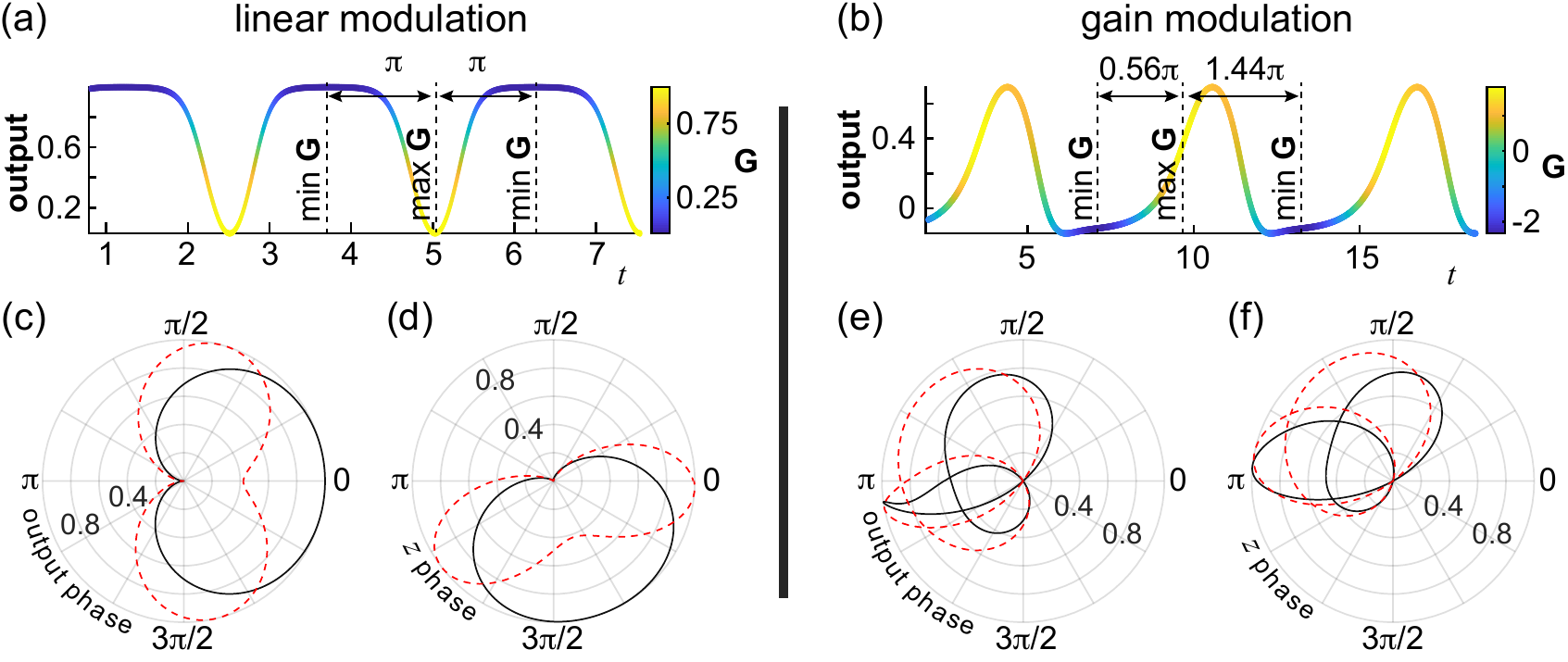}
\caption{\label{fig3}
Relationship between oscillation phase and instantaneous gain, $G$, in the two variations of the system. Minima/maxima of $G$ and their phase differences are annotated on the timeseries. The polar plots show the relationship between oscillation phase and normalized $G$ (divided by its maximum absolute value). Notice the constraining symmetry found in this relationship for the linearly modulated dynamics. Gain modulated dynamics exhibit a loose relationship and unconstrained phase plasticity, as indicated by the parameter perturbations (red dashed curves).
}
\end{figure*}

This interpretation of the qualitative difference can be visualized by the timeseries in Fig.~\ref{fig3}(a)-(b). The oscillation of the output is plotted across time, while instantaneous gain, $G$, is color-coded along the curve. Notice the one-to-one mapping between output and instantaneous gain in Fig.~\ref{fig3}(a): the same color appears across each value on the y-axis. Such constraint is not found in Fig.~\ref{fig3}(b), where gain modulation shapes the dynamics. Minima and maxima of $G$ are annotated and the oscillation phase difference is printed between consecutive extrema. Due to the fixed relationship between output and $G$, these phase differences equal $\pi$ in the case of linear modulation. In contrast, gain modulated dynamics exhibit no such constraint.

A deeper examination of the relationship between oscillation phase (extracted using Hilbert transform) and instantaneous gain is achieved through the polar plots in Fig.~\ref{fig3}. Instantaneous gain is plotted against the phase of output oscillation in Fig.~\ref{fig3}(c),(e) (black solid curves), corresponding to the timeseries in Fig.~\ref{fig3}(a)-(b). As expected, the relationship is symmetric around a specific phase value for the linearly modulated dynamics (Fig.~\ref{fig3}(c)). This is not the case for the oscillation produced through gain modulation (Fig.~\ref{fig3}(e), curve folded due to positive and negative values of $G$). Since the output oscillation is phase-locked to the oscillation of state variables, the same observations can be made by plotting $G$ against the phase of $z$'s, or any other state variable's, oscillation in Fig.~\ref{fig3}(d),(f).

Gain modulation, therefore, breaks the symmetry in the relationship between oscillation phase and instantaneous gain and this translates into enhanced phase plasticity and adaptability. This is demonstrated by the perturbed relationships indicated by the red dashed curves in Fig.~\ref{fig3}(c)-(f). These perturbations represent limit cycles produced after decreasing parameter $d$ by 10\% in both variations of the system. Notice the persistence of symmetry around the same phase value in Fig.~\ref{fig3}(c)-(d) and the elaborate, unconstrained, phase shift in Fig.~\ref{fig3}(e)-(f). Both variations of the system exhibit phase plasticity. However, it is constrained in the first case due to symmetry, while it is unleashed by gain modulation in the second.

\section*{Generalization}
\label{Gen}

The qualitative difference between the dynamics of the two variations, as seen in Figs.~\ref{fig2} and \ref{fig3}, cannot be attributed to anything else other than the different modulations shaping those dynamics. For instance, the different parameter sets used to produce the limit cycles cannot influence the symmetry found in the polar plots, as demonstrated by the perturbations in Fig.~\ref{fig3}(c)-(d). It is specifically the interactions between the variables and the distinct modulations of the I/O function that dictate the dimensionality of the dynamics in modulation phase space and, consequently, the relationship between oscillation phase and instantaneous gain. Furthermore, these properties are not limited to the dynamics of limit cycles but they apply to the whole repertoire of nonlinear dynamics in the system. It is exactly the same modulations that produce the fluctuations around a stable focus, for example.

The findings can be generalized to other systems by considering the modulations involved in the shaping of nonlinear dynamics in any given nonlinear system. The 2-dimensional dynamics observed in the modulation phase space of the first variation result from the shifting of the input-output curve along the input axis (see shading in Fig.~\ref{fig1}(c)). Whatever the position of the curve along that axis, it is obvious that the instantaneous gain (i.e. slope of the tangent) cannot change for a specific output value. Such linear modulations are generally described by $\bm{f_1}(x,z) = \bm{g}(z+\bm{h}(x))$, with modulator $\bm{h}(x)$ linearly modulating the \textit{input} of nonlinear function $\bm{g}$ (see below for \textit{output}'s linear modulation). Analytically, $G=\frac{\partial \bm{f_1}(x,z)}{\partial z} = \frac{\partial}{\partial z}\bm{g}(z+\bm{h}(x)) \times \frac{\partial}{\partial z}(z+\bm{h}(x)) = \frac{\partial}{\partial z}\bm{g}(z+\bm{h}(x))$, which is a function of $z+\bm{h}(x)$, as output $\bm{f_1}(x,z)$ is. If either $G$ or $\bm{f_1}(x,z)$ is a strictly monotonic function of $z+\bm{h}(x)$, or both, then this modulation yields an one-to-one mapping between output and instantaneous gain and, thus, 2-dimensional dynamics in modulation phase space.

\begin{figure*}
\includegraphics{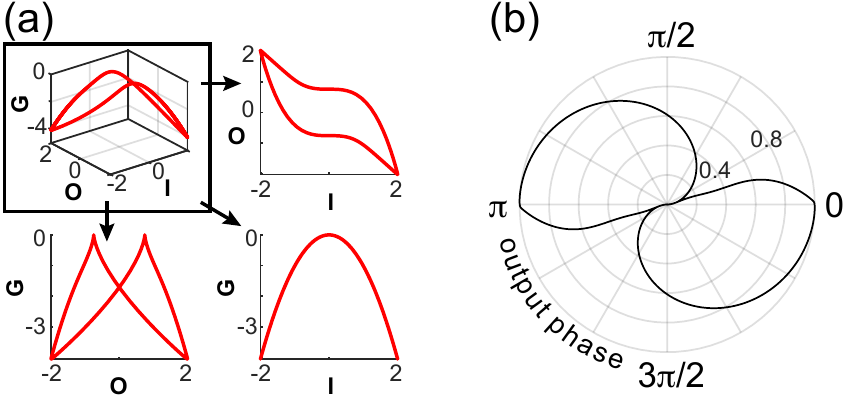}
\caption{\label{fig4}
Modulation phase space analysis of a limit cycle in the FitzHugh-Nagumo system for generalization purposes. The only multivariate nonlinear term, $-z^3/3 -x$, is analyzed as an example of linear modulation of the \textit{output} of a nonlinear I/O function (i.e. cubic in this case). (a) The limit cycle evolves on a 2-dimensional manifold in modulation phase space. Since it is produced through linear modulation of the \textit{output}, there is an one-to-one mapping between instantaneous gain and input. (b) The 2-dimensional constraint on the dynamics in modulation phase space translates into a symmetric relationship between instantaneous gain and oscillation phase of the \textit{output}, similar to Fig.~\ref{fig3}(c). 
}
\end{figure*}

Linear modulation can also be applied to the \textit{output} of a nonlinear function $\bm{g}$, described by $\bm{w}(x,z)=\bm{g}(z)+\bm{h}(x)$. This is the case for the FitzHugh–Nagumo system, for instance, where the output of the cubic I/O function of one variable is linearly modulated by another variable \cite{Fitzhugh1961}. Analytically, $G=\frac{\partial \bm{w}(x,z)}{\partial z} = \frac{\partial (\bm{g}(z)+\bm{h}(x))}{\partial z} = \frac{\mathrm{d}\bm{g}(z)}{\mathrm{d}z}$ which is a function of input $z$ for any nonlinear function $\bm{g}$. Therefore, the dynamics remain 2-dimensional in modulation phase space (Fig.~\ref{fig4}), however the instantaneous gain collapses onto input rather than the output, since the function is shifted along the output axis instead of the input axis. Consequently, the relationship between instantaneous gain and oscillation phase is constrained as in Fig.~\ref{fig3}(c)-(d).

\begin{figure*}
\includegraphics{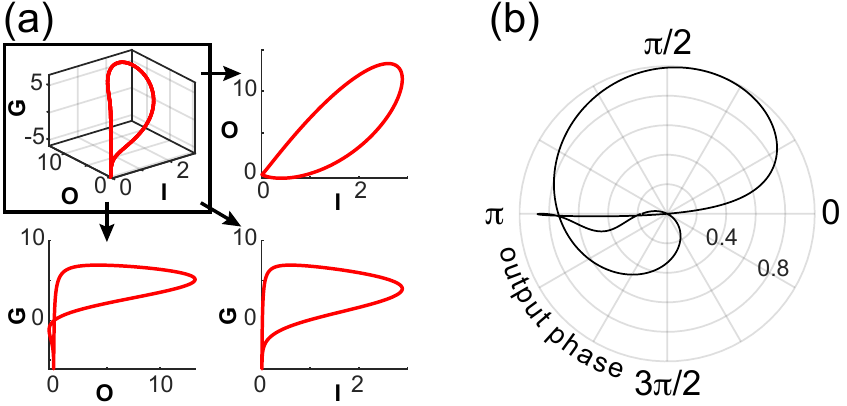}
\caption{\label{fig5}
Modulation phase space analysis of a limit cycle in the original R\"ossler system for generalization purposes. The only multivariate nonlinear term, $xz$, is analyzed as an example of gain modulation of the identity function, which is linear. (a) The limit cycle evolves on a 3-dimensional manifold in modulation phase space, similar to the second variation of the system. (b) The relationship between instantaneous gain and oscillation phase of the \textit{output} is free of any constraining symmetry, similar to Fig.~\ref{fig3}(e). 
}
\end{figure*}

Regarding the original R\"ossler system, with $\bm{f_0}(x,z) = xz$, it already features gain modulation applied on the \textit{default} linear I/O function, which is the identity function. That is, $x$ gain modulates $\bm{g_0}(z)=z$, or vice versa. Analyzing in respect to $z$, gives $G=x$, which is neither a function of input $z$ nor of output $\bm{f_0}(x,z)$. Thus, this modulation leads to 3-dimensional dynamics in modulation phase space (Fig.~\ref{fig5}), as in the second variation of the system. In this respect, the original R\"ossler system and its second variation are qualitatively similar. I purposefully imposed the extra nonlinearity of the sigmoidal to demonstrate both linear and gain modulations on the same input-output function, $\bm{g}$, without deteriorating the nonlinearity in the first variation of the system. That would be the case, since linear modulation of the identity function does not provide any nonlinearity to the system (nonlinear term $xz$ becomes linear $x+z$). The original R\"ossler system is therefore a minimal example of a nonlinear system featuring the simplest gain modulation possible: its only nonlinear term represents gain modulation of the identity function.

\section*{Discussion}

As demonstrated above, it is insightful to conceptualize the emergence of nonlinear dynamics through either linear modulations of nonlinear I/O functions, or gain modulations of even linear ones. Each type of modulation corresponds to a specific type of interactions between the variables in a nonlinear system. Linear modulations, which shift I/O functions along the input and output axes, produce 2-dimensional dynamics in modulation phase space (Fig.~\ref{fig2}), translating into a constrained relationship between instantaneous gain and the phase of dynamics (see symmetry in Fig.~\ref{fig3}). In contrast, gain modulations, which modulate the slope of I/O functions, produce 3-dimensional dynamics in modulation phase space (Fig.~\ref{fig2}) and loosen up the relationship between instantaneous gain and the phase of dynamics. The additional dimension in which dynamics are formed expands the adaptation potential of the system by unlocking phase plasticity, as demonstrated through parameter perturbations (Fig.~\ref{fig3}). Thus, gain modulation renders nonlinear dynamics highly malleable.

The universal applicability of the above conceptualization suggests that all multivariable nonlinear dynamical systems can be split into two classes. The first class includes all those systems that do not feature any gain modulation. The nonlinear dynamics are exclusively produced through linear modulation of nonlinear I/O functions, resulting into 2-dimensional dynamics in modulation phase space. Thus, the systems in the first class produce constrained nonlinear dynamics. The second class includes systems with at least one nonlinear term with gain modulation which can lead to 3-dimensional dynamics in modulation phase space. Thus, the systems in the second class produce loose nonlinear dynamics.

The present study adds to the literature dedicated to the entrainability and synchronizability of oscillatory dynamics. The findings indicate that gain modulation upgrades the adaptability of nonlinear dynamics, a property on which entrainment and synchronization relies on \cite{Perc2004}. The adaptability of oscillations and their synchronization depend on local contractive properties of the limit cycle, with low local divergence facilitating quick synchronization \cite{Perc2004, Marhl2006}. The state's vicinity to bifurcations has also been reported as a related factor \cite{Perc2004, Hasegawa2013}. Furthermore, the entrainability of oscillations is subserved by phase plasticity which indicates the adaptation of limit cycles in phase space \cite{Hatakeyama2015, Paijmans2017} (equivalently, phase compression subserving synchronizability \cite{Somers1993}). It is the same phase plasticity which is unlocked by gain modulation in this Letter. Thus, the present findings showcase gain modulation as a mechanism, not merely increasing adaptability, but fundamentally unfolding the adaptation possibilities in a nonlinear system. Consequently, this effect of gain modulation potentially explains the enhanced entrainment reported in a neural mass model, where gain modulation represented a neurophysiological mechanism \cite{Papasavvas2020}(see also \cite{Papasavvas2017} for its potential effect on neural synchronization).

Advances in nonlinear dynamics have the potential of informing all scientific disciplines in which nonlinear phenomena are integral. Understanding the role of gain modulation in such phenomena is beneficial. In fact, I would argue that such understanding is particularly important in neuroscience where gain modulation is a well documented neurophysiological mechanism \cite{Salinas2000, Ferguson2020}. Gain modulation is known to play a central role in neural computation and the processing of sensory information \cite{Carandini1994, Ayaz2009, Mitchell2003, Orton2016}. However, its impact on recurrent neural dynamics in cortical circuits has yet been sparsely investigated. Theoretical and experimental studies have started exploring the role of gain modulation in neural oscillations and their properties at cellular and circuit levels \cite{Shpiro2007, Kayser2015, Huguet2017, Papasavvas2015, Papasavvas2020}. The present findings suggest that the mechanism of gain modulation, delivered through divisive inhibition and other means \cite{Ferguson2020}, shapes neural oscillations in a way that allows for elaborate phase shift and high flexibility. Consequently, the oscillations' entrainment to rhythmic input and their synchronization across distributed brain networks are hypothesized to depend on the properties of the available gain modulation mechanisms, as applied to both cellular and circuit levels.

\section*{Acknowledgments}
I would like to thank Yujiang Wang, Peter Taylor, and Gerold Baier for discussion and comments on the manuscript. I am grateful to Andrew Trevelyan and Adrian Rees for introducing me to the concept of gain modulation in neuroscience. Thanks also to AP and NG for their generous hospitality during this study.

\section*{Supplementary Information}
Animations (GIF files) showing how the modulations of the sigmoidal I/O function form the limit cycles in each variation of the system (Fig.~\ref{fig1}). One shows limit cycle formation through a linear modulation and the other shows limit cycle formation through gain modulation.

\section*{Competing interests}
I declare no competing interests.



\bibliography{myBib}

\end{document}